%
%
%
%
%
%
\input amstex
\documentstyle {amsppt}
\topmatter
\title {{\it If there is an exactly $\lambda$-free abelian group then there is
an exactly $\lambda$-separable one in $\lambda$}\\
\medskip
Sh521} \endtitle
\rightheadtext {$\lambda$-free, $\lambda$-separable}
\author {Saharon Shelah \thanks {\null\newline
I thank Alice Leonhardt for the beautiful typing \newline
Typed December 12/93 - Done Fall '93;\S2 done 28 (+1) Feb.94 \newline
last revision 3/13/95} \endthanks} \endauthor
\affil {Institute of Mathematics \\
The Hebrew University \\
Jerusalem, Israel
\medskip
Rutgers University \\
Department of Mathematics \\
New Brunswick, NJ USA} \endaffil
\bigskip
\abstract {We give a solution stated in the title
to problem 3 of part 1 of the problems listed in the book of Eklof and
Mekler \cite{EM},(p.453).
There, in pp. 241-242, this is discussed and proved in some cases.
The existence of strongly $\lambda$-free ones was proved earlier by the
criteria in \cite{Sh:161} in \cite{MkSh:251}.
We can apply a similar proof to a large class of
other varieties in particular
to the variety of (non-commutative) groups.} \endabstract
\endtopmatter
\document

%
%
\def\renewcommand{\newcommand}	       
\edef\cite{\the\catcode`@}%
\catcode`@ = 11
\let\@oldatcatcode = \cite
\chardef\@letter = 11
\chardef\@other = 12
%
%
%
%
\def\@innerdef#1#2{\edef#1{\expandafter\noexpand\csname #2\endcsname}}%
%
%
\@innerdef\@innernewcount{newcount}%
\@innerdef\@innernewdimen{newdimen}%
\@innerdef\@innernewif{newif}%
\@innerdef\@innernewwrite{newwrite}%
%
%
%
\def\@gobble#1{}%
%
%
%
\ifx\inputlineno\@undefined
   \let\@linenumber = \empty 
\else
   \def\@linenumber{\the\inputlineno:\space}%
\fi
%
%
%
\def\@futurenonspacelet#1{\def\cs{#1}%
   \afterassignment\@stepone\let\@nexttoken=
}%
\begingroup 
\def\\{\global\let\@stoken= }%
\\ 
\endgroup
\def\@stepone{\expandafter\futurelet\cs\@steptwo}%
\def\@steptwo{\expandafter\ifx\cs\@stoken\let\@@next=\@stepthree
   \else\let\@@next=\@nexttoken\fi \@@next}%
\def\@stepthree{\afterassignment\@stepone\let\@@next= }%
%
%
%
\def\@getoptionalarg#1{%
   \let\@optionaltemp = #1%
   \let\@optionalnext = \relax
   \@futurenonspacelet\@optionalnext\@bracketcheck
}%
%
%
\def\@bracketcheck{%
   \ifx [\@optionalnext
      \expandafter\@@getoptionalarg
   \else
      \let\@optionalarg = \empty
      \expandafter\@optionaltemp
   \fi
}%
\def\@@getoptionalarg[#1]{%
   \def\@optionalarg{#1}%
   \@optionaltemp
}%
%
%
%
\def\@nnil{\@nil}%
\def\@fornoop#1\@@#2#3{}%
\def\@for#1:=#2\do#3{%
   \edef\@fortmp{#2}%
   \ifx\@fortmp\empty \else
      \expandafter\@forloop#2,\@nil,\@nil\@@#1{#3}%
   \fi
}%
\def\@forloop#1,#2,#3\@@#4#5{\def#4{#1}\ifx #4\@nnil \else
       #5\def#4{#2}\ifx #4\@nnil \else#5\@iforloop #3\@@#4{#5}\fi\fi
}%
\def\@iforloop#1,#2\@@#3#4{\def#3{#1}\ifx #3\@nnil
       \let\@nextwhile=\@fornoop \else
      #4\relax\let\@nextwhile=\@iforloop\fi\@nextwhile#2\@@#3{#4}%
}%
%
%
%
\@innernewif\if@fileexists
\def\@testfileexistence{\@getoptionalarg\@finishtestfileexistence}%
\def\@finishtestfileexistence#1{%
   \begingroup
      \def\extension{#1}%
      \immediate\openin0 =
         \ifx\@optionalarg\empty\jobname\else\@optionalarg\fi
         \ifx\extension\empty \else .#1\fi
         \space
      \ifeof 0
         \global\@fileexistsfalse
      \else
         \global\@fileexiststrue
      \fi
      \immediate\closein0
   \endgroup
}%
%
%
%
%
\def\bibliographystyle#1{%
   \@readauxfile
   \@writeaux{\string\bibstyle{#1}}%
}%
\let\bibstyle = \@gobble
%
%
\let\bblfilebasename = \jobname
\def\bibliography#1{%
   \@readauxfile
   \@writeaux{\string\bibdata{#1}}%
   \@testfileexistence[\bblfilebasename]{bbl}%
   \if@fileexists
      \nobreak
      \@readbblfile
   \fi
}%
\let\bibdata = \@gobble
%
%
\def\nocite#1{%
   \@readauxfile
   \@writeaux{\string\citation{#1}}%
}%
\@innernewif\if@notfirstcitation
%
%
\def\cite{\@getoptionalarg\@cite}%
%
%
\def\@cite#1{%
   \let\@citenotetext = \@optionalarg
   \printcitestart
   \nocite{#1}%
   \@notfirstcitationfalse
   \@for \@citation :=#1\do
   {%
      \expandafter\@onecitation\@citation\@@
   }%
   \ifx\empty\@citenotetext\else
      \printcitenote{\@citenotetext}%
   \fi
   \printcitefinish
}%
\def\@onecitation#1\@@{%
   \if@notfirstcitation
      \printbetweencitations
   \fi
   \expandafter \ifx \csname\@citelabel{#1}\endcsname \relax
      \if@citewarning
         \message{\@linenumber Undefined citation `#1'.}%
      \fi
      \expandafter\gdef\csname\@citelabel{#1}\endcsname{%
         {\tt
            \escapechar = -1
            \nobreak\hskip0pt
            \expandafter\string\csname#1\endcsname
            \nobreak\hskip0pt
         }%
      }%
   \fi
   \csname\@citelabel{#1}\endcsname
   \@notfirstcitationtrue
}%
%
%
\def\@citelabel#1{b@#1}%
%
%
\def\@citedef#1#2{\expandafter\gdef\csname\@citelabel{#1}\endcsname{#2}}%
%
%
%
\def\@readbblfile{%
   \ifx\@itemnum\@undefined
      \@innernewcount\@itemnum
   \fi
   \begingroup
      \def\begin##1##2{%
         \setbox0 = \hbox{\biblabelcontents{##2}}%
         \biblabelwidth = \wd0
      }%
      \def\end##1{}
      %
      %
      \@itemnum = 0
      \def\bibitem{\@getoptionalarg\@bibitem}%
      \def\@bibitem{%
         \ifx\@optionalarg\empty
            \expandafter\@numberedbibitem
         \else
            \expandafter\@alphabibitem
         \fi
      }%
      \def\@alphabibitem##1{%
         \expandafter \xdef\csname\@citelabel{##1}\endcsname {\@optionalarg}%
         \ifx\biblabelprecontents\@undefined
            \let\biblabelprecontents = \relax
         \fi
         \ifx\biblabelpostcontents\@undefined
            \let\biblabelpostcontents = \hss
         \fi
         \@finishbibitem{##1}%
      }%
      \def\@numberedbibitem##1{%
         \advance\@itemnum by 1
         \expandafter \xdef\csname\@citelabel{##1}\endcsname{\number\@itemnum}%
         \ifx\biblabelprecontents\@undefined
            \let\biblabelprecontents = \hss
         \fi
         \ifx\biblabelpostcontents\@undefined
            \let\biblabelpostcontents = \relax
         \fi
         \@finishbibitem{##1}%
      }%
      \def\@finishbibitem##1{%
         \biblabelprint{\csname\@citelabel{##1}\endcsname}%
         \@writeaux{\string\@citedef{##1}{\csname\@citelabel{##1}\endcsname}}%
         \ignorespaces
      }%
      %
      %
      \let\em = \bblem
      \let\newblock = \bblnewblock
      \let\sc = \bblsc
      \frenchspacing
      \clubpenalty = 4000 \widowpenalty = 4000
      \tolerance = 10000 \hfuzz = .5pt
      \everypar = {\hangindent = \biblabelwidth
                      \advance\hangindent by \biblabelextraspace}%
      \bblrm
      \parskip = 1.5ex plus .5ex minus .5ex
      \biblabelextraspace = .5em
      \bblhook
      \input \bblfilebasename.bbl
    \endgroup
}%
%
%
\@innernewdimen\biblabelwidth
\@innernewdimen\biblabelextraspace
%
%
%
\def\biblabelprint#1{%
   \noindent
   \hbox to \biblabelwidth{%
      \biblabelprecontents
      \biblabelcontents{#1}%
      \biblabelpostcontents
   }%
   \kern\biblabelextraspace
}%
%
%
%
\def\biblabelcontents#1{{\bblrm [#1]}}%
%
%
\def\bblrm{\rm}%
%
%
\def\bblem{\it}%
%
%
\def\bblsc{\ifx\@scfont\@undefined
              \font\@scfont = cmcsc10
           \fi
           \@scfont
}%
%
%
\def\bblnewblock{\hskip .11em plus .33em minus .07em }%
%
%
\let\bblhook = \empty
%
%
%
\def\printcitestart{[}
\def\printcitefinish{]}
\def\printbetweencitations{, }
\def\printcitenote#1{, #1}
%
%
%
\let\citation = \@gobble
%
%
%
\@innernewcount\@numparams
%
%
\def\newcommand#1{%
   \def\@commandname{#1}%
   \@getoptionalarg\@continuenewcommand
}%
%
%
\def\@continuenewcommand{%
   \@numparams = \ifx\@optionalarg\empty 0\else\@optionalarg \fi \relax
   \@newcommand
}%
%
%
\def\@newcommand#1{%
   \def\@startdef{\expandafter\edef\@commandname}%
   \ifnum\@numparams=0
      \let\@paramdef = \empty
   \else
      \ifnum\@numparams>9
         \errmessage{\the\@numparams\space is too many parameters}%
      \else
         \ifnum\@numparams<0
            \errmessage{\the\@numparams\space is too few parameters}%
         \else
            \edef\@paramdef{%
               \ifcase\@numparams
                  \empty  No arguments.
               \or ####1%
               \or ####1####2%
               \or ####1####2####3%
               \or ####1####2####3####4%
               \or ####1####2####3####4####5%
               \or ####1####2####3####4####5####6%
               \or ####1####2####3####4####5####6####7%
               \or ####1####2####3####4####5####6####7####8%
               \or ####1####2####3####4####5####6####7####8####9%
               \fi
            }%
         \fi
      \fi
   \fi
   \expandafter\@startdef\@paramdef{#1}%
}%
%
%
%
%
\def\@readauxfile{%
   \if@auxfiledone \else 
      \global\@auxfiledonetrue
      \@testfileexistence{aux}%
      \if@fileexists
         \begingroup
            \endlinechar = -1
            \catcode`@ = 11
            \input \jobname.aux
         \endgroup
      \else
         \message{\@undefinedmessage}%
         \global\@citewarningfalse
      \fi
      \immediate\openout\@auxfile = \jobname.aux
   \fi
}%
%
%
\newif\if@auxfiledone
\ifx\noauxfile\@undefined \else \@auxfiledonetrue\fi
%
%
%
%
\@innernewwrite\@auxfile
\def\@writeaux#1{\ifx\noauxfile\@undefined \write\@auxfile{#1}\fi}%
%
%
%
\ifx\@undefinedmessage\@undefined
   \def\@undefinedmessage{No .aux file; I won't give you warnings about
                          undefined citations.}%
\fi
%
%
\@innernewif\if@citewarning
\ifx\noauxfile\@undefined \@citewarningtrue\fi
%
%
%
\catcode`@ = \@oldatcatcode

%
%
\newpage

\head {\S0 Introduction} \endhead
\bigskip

\subhead {Convention} \endsubhead  In \S0 and \S1, ``group" here means
``abelian group", and ``free" means in this variety. \newline

We assume there is a $\lambda$-free, non-free (abelian) group of cardinality
$\lambda$.  We shall prove that there is a $\lambda$-separable non-free
abelian group of cardinality $\lambda$, apriori a stronger statement.  We
rely on the characterization of $\lambda$ as in the hypothesis from 
\cite{Sh:161}:
the existence of $S,\langle <s^\ell_\eta:\ell < n>:\eta \in S_f \rangle$,
$\langle \ell(k):k \rangle$
as there (see appendix; i.e. \S3 here).  Mekler Shelah \cite{MkSh:251} dealt
with a similar weaker problem in a parallel way: if there is a 
$\lambda$-free not free
abelian group of cardinality $\lambda$ then there is a strongly 
$\lambda$-free one.  In Eklof Mekler \cite{EM}, the present problem 
was raised, discussed and sufficient conditions were given, depending 
on the form of $S$, see \cite{EM}, p.242-242, the problem in \cite{EM},p.453.
The direct sufficient condition is that for every $S' \subseteq S_f$ of
cardinality $\lambda$ there is a well ordering $<^*$ such that for each
$\eta \in S_f, \dsize \bigcup_{\ell < n} s^\ell_\eta$ is almost disjoint to
$\cup\{ \dsize \bigcup_{\ell < n} s^\ell_\nu:\nu <^* \eta$ and $\nu \in S\}$.
In particular from the assumption for $\lambda$,
the conclusion for $\lambda^+$ (i.e. the existence of such $S$) was gotten.
However, not all cases were done there.  Our approach is more
algebraic.  In \S2 we deal with generalizations to other varieties and in
\S3 we present relevant material from \cite{Sh:161} (on $\lambda$-systems)
to make the paper self-contained.
\bigskip

\subhead {Explanation of the proof of the main theorem} \endsubhead  
It may be helpful to read this explanation if you are lost or stuck during
the proof but it assumes some notations from the proof.
We construct $G$ that is freely generated by
$x[a]$ (for $a \in \dsize \bigcup_{\eta \in S_c} B_\eta)$ and
$y_{\eta,m}$ (for $\eta \in S_f$ and $m < \omega$) except the equation

$$
2y_{\eta,m+1} = y_{\eta,m} + \Sigma\{x[a^\ell_{\eta,m}]:\ell < n \}.
\tag "{$(*)_{\eta,m}$}"
$$
\medskip

\noindent
Let $G = G_{I_0}$, $I_0 = I_{<>,\lambda}$.

Let $\alpha < \lambda$ and we want to show that if $\alpha < \lambda$ and
$\langle \alpha \rangle \notin S$ then $G_{<>,\alpha}$ (which is essentially 
the subgroup
generated by the $y_{\eta,m}$ and $x[a^\ell_{\eta,m}]$ satisfying $\eta(0) 
< \alpha)$ is a free direct summand of $G = G_{<>,\lambda}$.

\noindent
We do not see combinatorially why this holds, so we find $I_1 \supseteq
 I_{<>,\alpha}$, $I_1 \in K^+$ such that

$$
\eta \in S_f \backslash S_f[I_1] \Rightarrow
\dsize \bigcup_{\ell < n} s^\ell_\eta \text{ is almost disjoint to }
Y[I_1] \tag{$**$}
$$
\medskip

\noindent
So let
$g_{I_0,I_1}$ be the natural homomorphism from $G_{I_0}$ to
$G_{I_1}$; well, why does it work? by $(**)$.

\noindent
Also $g_{I_0,I_1}$ is the identity on $G_{<>,\alpha}$ and $G_{I_1} \backslash
G_{<>,\alpha}$ is $\cong G_{I_2}$ where $I_2 = I_1 \backslash
I_{<>,\alpha}$,
but $I_2 \in K^+$ so $G_{I_1}/G_{<>,\alpha}$ is free hence $G_{<>,\alpha}$ is
a direct summand of $G_{I_1}$, so there is a 
projection $f$ from $G_{I_1}$ onto
$G_{I_{\langle \rangle,\alpha}}$ so $f \circ g_{I_0,\alpha}$ is a projection
from $G$ onto $G_{I_{\langle \rangle,\alpha}}$ and we can
complete the proof.

To accomplish $(**)$ we need good control over how e.g. $s^\ell_\eta$ \, 
$(\eta(0) > \alpha)$ intersect $B_{<\alpha>}$, and this
is the information we put in the appendix on the $\lambda$-system 
(really old \cite{Sh:161} is O.K., but we retain the appendix to ease 
reading).
\bigskip

\demo{0.1 Definition}  For $\Xi$ a set of variables, $\Gamma$ set of
equations in some variables (maybe outside $\Xi$)
let $G(\Xi,\Gamma)$ be the (abelian) group
freely generated by $\Xi$, except the equations in $\Gamma \restriction \Xi$,
i.e. the equations from $\Gamma$ mentioning only variables from $\Xi$.
\enddemo
\bigskip

\demo{0.2 Observation}  1) A sufficient condition (assuming $\Xi \subseteq
\Xi'$ sets of variables) for
\medskip
\roster
\item "{$(*)$}"  $G(\Xi',\Gamma)$ is a free extension of $G(\Xi,\Gamma)$ (i.e.
the mapping induced by $id_\Xi$ from $G(\Xi,\Gamma)$ into $G(\Xi',\Gamma)$,
which is always homomorphism, is an embedding, and $G(\Xi^\prime,\Gamma)$ 
divided by the range of this mapping is a free group),
\endroster
\medskip

\noindent
is
\medskip
\roster
\item "{$(**)$}"  there is an increasing continuous sequence
 $\langle \Xi_\zeta:\zeta \le \zeta^* \rangle$, $\Xi_0 = \Xi,\Xi_{\zeta^*} =
\Xi'$, and $G(\Xi_{\zeta + 1},\Gamma)$ is a free extension of $G(\Xi_\zeta,
\Gamma)$.
\endroster
\medskip

\noindent
2)  Another sufficient condition for $(*)$ of 0.2, is that by change of the
variables in $\Xi' \backslash \Xi$, the set of equations $\Gamma 
\restriction \Xi'$ is only $\Gamma \restriction \Xi$.
\enddemo
\newpage

\head {\S1 Proving $\lambda$-Separability} \endhead
\bigskip

Here we prove the main theorem; the reader is advised to look at 3.6, 3.7
at least during reading the beginning of the proof, and also to look again
at the explanation in \S0 of the proof when arriving to read the middle of
the proof.
\bigskip

\definition {1.1 Definition} A group $G$ is $\lambda$-separable if:\newline
$H \subseteq G,Rk(H) < \lambda \Rightarrow H$ included in a free direct
summand of $G$. \newline
(Remember: for an uncountable group $H$, its rank, $Rk(H)$ is equal to its
cardinality, $|H|$.)
\enddefinition
\bigskip

\proclaim {1.2 Main Theorem}  If there is a $\lambda$-free non
$\lambda^+$-free
(abelian) group $(\lambda > \aleph_0,\lambda$ necessarily regular)
\underbar{then} there
is a $\lambda$-free, $\lambda$-separable, not $\lambda^+$-free group.
\endproclaim
\bigskip

\demo {Proof}  The hypothesis of the theorem on the existence of such
groups is analyzed in detail in \cite{Sh:161} (most relevant are
\cite{Sh:161},3.6,3.7), and in particular, it implies the existence of
$n,S,\lambda(\eta,S),\langle B_\eta:\eta \in S_c \rangle$,
$\langle s^\ell_\eta:\eta \in S_f,\ell < n \rangle$ with the properties as in
\cite{Sh:161},3.6,3.7 presented in 3.6, 3.7 of the Appendix here,
and let \newline
$\langle a^\ell_{\eta,m}:m < \omega \rangle$ list $s^\ell_\eta$ in
increasing order for the order
of $B_{(\eta \restriction \ell) \char 94 \langle \lambda(\eta
\restriction \ell) \rangle}$ (see clause (i) of 3.7) and without loss 
of generality we have in addition
\medskip
\roster
\item "{$(*)$}"  for $\eta \in S_f,\ell < n$, we have \newline
$\alpha_{\eta,\ell,m} =: \text{ min}\{ \beta:a^\ell_{\eta,m} \in 
B_{(\eta \restriction \ell) \char 94 <\beta>} \} < \lambda(\eta \restriction
\ell,S)$ \newline
is non-decreasing in $m$,
\endroster
and we call its limit $\beta^*(\eta,\ell)$ (so
$s^\ell_\eta \subseteq B_{(\eta \restriction \ell) \char 94 \langle \beta^*
(\eta,\ell) \rangle}$ and $\beta^*(\eta,\ell) \le \eta(\ell)$).
\medskip
\roster
\item "{$(**)$}"  if $\rho \in S_i,\nu \triangleleft \rho,k = \ell g(\nu)$
and $\text{cf}(\rho(k)) = \lambda(\rho,S)$ then
{\roster
\itemitem{ (a) }  for $\beta < \lambda(\rho,S)$ we have
$\text{sup}\{ \beta^*(\eta,k):
\rho \char 94 \langle \beta \rangle \trianglelefteq \eta \in S_f \}$ is
$< \rho(k)$
\itemitem{ (b) }  the sequence $\langle \text{min} \{ \beta^*(\eta,k):
\rho \char 94 \langle \beta \rangle \trianglelefteq \eta \in S_f \}:\beta <
\lambda(\rho,S) \rangle$ is strictly increasing with limit $\rho(k)$.
\endroster}
\endroster
\medskip

\noindent
(see Appendix, clauses $(f)(\alpha),(f)(\beta)$ and (g) of 3.6). \newline
Let
$$
\align
K = \{I:&I \subseteq S_c \text{ and } 
[\eta \ne \nu \and \eta \in I \and \nu \in
I \Rightarrow \neg(\eta \trianglelefteq \nu)] \text{ and } \\
  &[\eta \char 94 <\beta> \, \in I \and \alpha < \beta 
\Rightarrow \eta \char 94 <\alpha> \, \in I]\},
\endalign
$$

$$
K^+ = \{ I \in K:I \ne \emptyset \text{ and } \eta \char 94 <\alpha> \in I
\Rightarrow \vee_{\beta < \lambda(\eta,I)}[\eta \char 94 <\beta> \,
\notin I] \}.
$$

\noindent
For $I \in K$ let $J[I] =: \{ \eta \in S_i$: for some $\alpha,\eta \char 94
<\alpha> \, \in I \}$, so for $\eta \in J[I]$ there is a unique $\alpha_I
[\eta] \le \lambda (\eta,S)$ such that $[\eta \char 94 <\alpha> \, \in I
\Leftrightarrow \alpha < \alpha _I[\eta]]$, note:\newline
$I \in K^+,\eta \in J[I]
\Rightarrow \alpha_I[\eta] < \lambda(\eta,S)$.  For $I \in K$ let \newline
$S_f[I] =: \{ \eta \in S_f:\text{for some } k, \eta \restriction k \in I \}$
; note that the $k$ is unique and if \newline
$I \ne \{<>\}$, then $k > 0$, so we choose to write
$k = k_I(\eta) + 1$ (so for \newline
$I = \{<>\},k_I(\eta) = -1)$.  Also let
$$
Y[I] =: \cup \{ B_\nu:\text{for some } \eta \in I \text{ we have }
\eta \trianglelefteq \nu \in S_c \}
$$

\noindent
For $\eta \in S_f$ let $w_I(\eta) = \{ \ell < n:\text{for every }
m < \omega,a^\ell_{\eta,m} \in Y[I] \}$ (equivalently: for infinitely
many $m < \omega,a^\ell_{\eta,m} \in Y[I]$). \newline
For every $I \in K$ we define a group $G_I$, it is freely generated by
\newline
$\Xi_I =: \{y_{\eta,m}:m < \omega$ and $\eta \in S_f[I]\} \cup \{x[a]:
a \in Y[I]\}$ except the equations (we call this set $\Gamma_I$):
\medskip
\roster
\item "{$(*)_I$}"  $\quad$ for $\eta \in S_f[I]$ and $m < \omega$, the 
equation $\varphi^m_{I,\eta}$ defined as:
\item "{$(*)^m_{I,\eta}$}"  $\quad 2y_{\eta,m+1} = y_{\eta,m} + \sum \{
x[a^\ell_{\eta,m}]:\ell < n$ and $a^\ell_{\eta,m} \in Y[I] \}$
\endroster
\medskip

\noindent
Note that ${}^1 \lambda \in K$, and let $G = G_{({}^1 \lambda)}$; this abelian
group is the example as in \newline
\cite{Sh:161}, Lemma 5.3, in particular $G$ is not free.
Let $<_{\ell x}$ be lexicographic order of $S$, clearly
it is a well ordering.
\enddemo
\bigskip

\subhead{A Fact} \endsubhead  If $I \in K^+$ then $G_I$ is free.
\bigskip

\demo{Proof}  We can find functions $\ell$ and $m$, where for
$\eta \in S_f[I]$ we have \newline
$\ell(\eta) \in \{ k_I(\eta),\dotsc,n-1 \}$ and $m(\eta) < \omega$
and we can find a list $\langle \eta_\zeta:\zeta < \zeta^* \rangle$
of $S_f[I]$ such that:
\medskip
\roster
\item "{$(*)$}"  $\{ a^{\ell(\eta_\zeta)}_{\eta_\zeta,m}:m \in [m(\eta_\zeta),
\omega)\}$ is disjoint to $\{a^\ell_{\eta_\varepsilon,m}:
a^\ell_{\eta_\varepsilon,m} \in Y[I],m < \omega$, \newline
$\varepsilon < \zeta,\ell < n \}$ \underbar{and} $\{ \eta \in S_f[I]:\eta
<_{\ell x} \nu \}$ is an initial segment of $\langle \eta_\zeta:\zeta <
\zeta^* \rangle$ for each $\nu \in J[I]$
\endroster
\medskip

\noindent
[why? for each $\nu \in J[I]$ well order $\{ \eta \in S_f[I]:\nu \triangleleft
\eta \}$ by \cite{Sh:161}, 3.10 (and 3.6 clause $(c)$ and the definition of 
$K^+$), say by $<^*_\nu$, then order the blocks by $<_{\ell x}$].
Without loss of generality $m(\eta)$ is minimal such that $(*)$ holds.

For $\zeta \le \zeta^*$ let $H_\zeta$ be the subgroup of $G_I$ generated 
by \newline
$\Xi_\zeta = \{ x[a^\ell_{\eta_\varepsilon,m}]:\varepsilon
 < \zeta,m < \omega,\ell \in [k(\eta),n)\} \cup \{ y_{\eta_\varepsilon,m}:
\varepsilon < \zeta,m < \omega \}$.  Let $H_{\zeta^* + 1} = G_I$.

Now $\langle H_\zeta:\zeta \le \zeta^* + 1 \rangle$ is increasing continuous,
$H_0 = \{ 0 \}$, $H_{\zeta^* + 1} = G_I$ and $H_{\zeta + 1}/H_\zeta$ is free.
Why?  we use 0.2(1), so
it is enough to prove $G(\Xi_{\zeta + 1},\Gamma_I)$
is a free extension of $G(\Xi_\zeta,\Gamma_I)$ for each $\zeta \le \zeta^*$.
For $\zeta = \zeta^*$, we just
add variables $(\{ x[a]:a \in \Xi_I \backslash \Xi_{\zeta^*} \})$ but no
equations.  For $\zeta < \zeta^*$, we can "forget" $y_{\eta_\zeta,m}$ for
$m < m(\eta_\zeta)$ and replace/omit $x[a^{\ell(\eta_\zeta)}_{\eta_\zeta,m}]$
for $m \in [m(\eta_\zeta),\omega)$,
so $G(\Xi_{\zeta + 1},\Gamma_I)$ is freely generated over
$G(\Xi_\zeta,\Gamma_I)$ by \newline
$\{ y_{\eta_\zeta,m}:m \ge m(\eta_\zeta)\} \cup 
\{ x[a]:x[a] \in \Xi_{\zeta + 1} \backslash \Xi_\zeta \backslash
\{ x[a^{\ell(\eta_\zeta)}_{\eta_\zeta,m}]:m \ge m(\eta_\zeta) \}\}$.
\enddemo
\bigskip

\subhead {B Notation} \endsubhead  Let $I_{\eta,\alpha} =: \{ \eta \char 94
<\beta>:\beta < \alpha \}$ (for $\eta \in S_i$).  Let $G_{\eta,\alpha} =:
G_{I_{\eta,\alpha}}$ so:
\medskip
\roster
\item "{{}}"  $G_{<>,\lambda}$ is the group $G$ (which we shall prove 
exemplifies the \newline

$\qquad$ conclusion of 1.2) \newline
$J[I_{<>,\lambda}] = \{ <> \}$
\endroster
\bigskip

\definition{C Definition}
1) $I_1 \le I_2$ (from $K$) if $S_f[I_1] \subseteq S_f
[I_2]$ and $(\forall \eta \in S_f(I_1))[k_{I_1}(\eta) \ge k_{I_2}(\eta)]$.
This implies $Y[I_1] \subseteq Y[I_2]$ and there is a clear relation between
$\Gamma_{I_1}$ and $\Gamma_{I_2}:\,$ each equation $\varphi = \varphi^m
_{I_1,\eta}$ in 
$\Gamma_{I_1}$ ``appears" in $\Gamma_{I_2}$ as $\psi = \varphi^m_{I_2,\eta}$
\underbar{but} $\psi$
is with more $x[a^\ell_{\eta,m}]$'s (for same old $\eta$
but new $\ell$'s which appear ``because" of some $\nu \in S_f[I_2] \backslash
S_f[I_1]$) and $\Gamma_{I_2}$ has members 
(not related to any equation from $\Gamma_{I_1}$) involving a new $\eta$.  
Another way to state this relation is
\newline
$(\forall \eta \in I_1)(\exists \nu \in I_2)[\nu \trianglelefteq \eta]$.
\newline
2)  $I_1 \le^d I_2$ \underbar{if} $I_1 \le I_2$ and 
$J[I_1]$ is a $<_{\ell x}$-initial segment of $J[I_2]$. 
\enddefinition
\bigskip

\subhead{D Fact}  \endsubhead  1)  $\le$ and $\le^d$ are partial orders
(of $K$). \newline
2)  If $I \in K \backslash \{\{ <> \}\}$ then $I = \dsize \bigcup
_{\eta \in J[I]} I_{\eta,\alpha_I(\eta)}$. \newline
3)  If $I_1 \le^d J_2$ then $Y[I_1]$ is a subset of $Y[I_2]$.
\bigskip

\demo{Proof}  Check.
\enddemo
\bigskip

\definition{E Definition}  Assume $I_1 \le^d I_2$ (both in $K$),
let $h_{I_1,I_2}$ be the
homomorphism from $G_{I_1}$ into $G_{I_2}$ defined by $h(x[a]) = x[a],
h(y_{\eta,m}) = y_{\eta,m}$ for $x[a] \in Y[I_1],\eta \in S_f[I_1]$.
\enddefinition
\bigskip

\definition{F Fact}  $h_{I_1,I_2}$ is really a homomorphism.
\enddefinition
\bigskip

\demo{Proof}  Look at the relevant equations.
\enddemo
\bigskip

\subhead{G Fact}  \endsubhead  If $I_1 \le^d I_2$ are from $K^+$ and
$(\forall \eta) \left[ \eta \in J[I_1] \Rightarrow \eta
\char 94 \langle \alpha_{I_1}(\eta) \rangle \notin S \right]$ \underbar{then}
\medskip
\roster
\item "{$(\alpha)$}"  $G_{I_2}/h_{I_1,I_2}(G_{I_1})$ is free and
\item "{$(\beta)$}"  $h_{I_1,I_2}$ is one to one.
\item "{$(\gamma)$}"  $\text{Rang}(h_{I_1,I_2}) = \langle y_{\eta,m},
x[a]:a \in Y[I_1],\eta \in S_f[I_1]$ and $m < \omega \rangle_{G_{I_2}}$
(i.e. the subgroup generated by this set)
\endroster
\medskip

\noindent
so we look at $h_{I_1,I_2}$ as the identity.
\bigskip

\demo{Proof}  Like the proof of Fact A.
\enddemo
\bigskip

\definition{H Conclusion}  If $I_1 \le^d I_2$
(so are from $K$) \underbar{then} $h_{I_1,I_2}$ is an embedding.
\enddefinition
\bigskip

\demo{Proof}  As a direct limit of ones satisfying the assumptions of Fact G.
\enddemo
\bigskip

\subhead{I Fact} \endsubhead
\medskip
\roster
\item "{$(\alpha)$}"  $G = G_{({}^1 \lambda)} = G_{<>,\lambda} = \dsize 
\bigcup_{\alpha < \lambda}G_{<>,\alpha}$ (increasing continuous)
\item "{$(\beta)$}"  for $\alpha < \lambda$ the group $G_{<>,\alpha}$ is free.
\endroster
\bigskip

\demo{Proof}  For clause $(\alpha)$ as $\Gamma_{({}^1 \lambda)} = 
\dsize \bigcup_{\alpha < \lambda} \left( \Gamma
\restriction \{ y_{\eta,m},x[a^\ell_{\eta,m}]:\eta(0) < \alpha \} \right)$,
using Fact H (see Fact $G$ last line). \newline
For clause $(\beta)$ see Fact $A$.
\enddemo
\bigskip

\definition{J Definition}  For $I_1 \le I_2$ (in $K$), satisfying
$\otimes_{I_1,I_2}$ below, \underbar{let}
$g_{I_2,I_1}$ be the homomorphism from $G_{I_2}$ into $G_{I_1}$ defined by:
\medskip
\roster
\item "{(i)}"  if $a \in Y[I_1]$ \underbar{then} $g_{I_2,I_1}(x[a]) = x[a]$
\item "{(ii)}"  if $a \in Y[I_2] \backslash Y[I_1]$ \underbar{then}
 $g_{I_2,I_1}(x[a]) = 0$
\item "{(iii)}"  if $\eta \in S_f[I_1]$ \underbar{then}
$g_{I_2,I_1}(y_{\eta,m}) = y_{\eta,m}$
\item "{(iv)}"  if $\eta \in S_f[I_2] \backslash S_f[I_1]$ and
$\{ a^\ell_{\eta,k}:k \in [m,\omega)$ and $\ell \ge k_{I_2}[\eta]$
(equivalently $a^\ell_{\eta,k} \in Y[I_2])\}$
is disjoint to $Y[I_1]$ \underbar{then} $g_{I_2,I_1}(y_{\eta,m}) = 0$
\endroster
\medskip

\noindent
(this is enough for defining $g_{I_2,I_1}$) \newline
where
\medskip
\roster
\item "{$\otimes_{I_1,I_2}$}"  for $\eta \in S_f[I_2] \backslash S_f[I_1]$,
$\dsize \bigcup_{\ell \in [k_{I_2}(\eta),n)} s^\ell_\eta$ is almost 
disjoint to $Y[I_1]$ \newline
(i.e. has finite intersection).
\endroster
\enddefinition
\bigskip

\subhead{K Fact}  \endsubhead  Assume $I_1 \le I_2$ are in $K$.  Then  
\roster
\item "{$(\alpha)$}"  $g_{I_2,I_1}$ really defines a homomorphism which is
onto (when $I_1 \le I_2$ and $\otimes_{I_1,I_2}$ holds)
\item "{$(\beta)$}"  Kernel$(g_{I_2,I_1})$ is the subgroup of $G_{I_2}$
generated by the set of $x[a]$'s and $y_{\eta,m}$'s which by Definition J
are sent by $g_{I_2,I_1}$ to $0$.
\endroster
\bigskip

\demo{Proof}  Check the equations.
\enddemo
\bigskip

\subhead {Main Fact L} \endsubhead  If $\alpha < \lambda$ and
$<\alpha> \, \notin S$
then $G_{<>,\alpha}$ is a direct summand of $G = G_{\{ <> \}}$.
\bigskip

\demo{Proof}  We can define by induction on $k$ a number
$\ell_k \le n:\ell_k = 0$,
if $\ell_k$ is defined and $< n$, let $\ell_{k+1}$ be the unique $\ell$
such that $\ell_k < \ell \le n$ and $\eta \in S_f \Rightarrow 
cf(\eta(\ell_k)) = \lambda(\eta \restriction \ell,
S)$ (exists by 3.3(f), all $\eta \in S_f$ behave the same by
3.6(a) (and see 3.2(6)(d)), note: if $\eta \in S_f \Rightarrow \text{ cf}
(\eta(\ell_k)) = \aleph_0$ then $\ell_{k+1} = n$.  Clearly if $\ell_k$ is
defined and $< n$ then $\ell_k < \ell_{k+1} \le n$.   So for some
$k^*$, $\ell_{k^*} = n$. \newline
We shall define by induction on $k \le k^*$ the following
$J_k$ and, when $k < k^*$, \newline
$\langle \alpha_\eta:\eta \in J_k \rangle$
such that:
\medskip
\roster
\item "{(0)}"  $J_k \subseteq S \cap {}^{\ell_k} \lambda$
\item "{(1)}"  $\alpha_\eta < \lambda(\eta,S)$ and $\eta \char 94 \langle
 \alpha_\eta \rangle \notin S$ for $\eta \in J_k$
\item "{(2)}"  $J_{k+1} = \{ \eta:\eta \in S \cap {}^{\ell_{k+1}}\lambda$ and
$\eta \restriction \ell_k \in J_k$ but $\eta(\ell_k) > 
\alpha_{\eta \restriction \ell_k} \}$
\item "{(3)}"  if $\eta \in J_{k+1},k+1 < k(*),\alpha \in [\alpha_\eta,
\lambda(\eta,S))$ and $\eta \char 94 <\alpha> \trianglelefteq \,\, 
\nu \in S_f$
\underbar{then} $s^{\ell_k}_\nu \cap B_{\eta \restriction \ell_k \char 94
\langle \alpha_{\eta \restriction \ell_k} \rangle}$ is finite
\item "{(4)}"  $J_0 = \{<> \}, \alpha_{<>} = \alpha$.
\endroster
\medskip

For $k=0$ use clause (4).  For $k+1$ we define $J_{k+1}$ by clause (2), now
if $k+1 < k(*)$ for $\eta \in J_{k+1}$ we have to find $\alpha_\eta$ to
satisfy clauses (1), (3), this is possible by $(*)$,$(**)$ in the beginning of
the proof of Theorem 1.2. \newline

Let $I_0 = \{ <\beta>:\beta < \lambda \}$,\newline
$I_1 = \{ \eta \char 94 <\beta>:\text{for some } k < k^* \text{ we have }
\eta \in J_k
\text{ and } \beta < \alpha_\eta \}$, \newline
$I_2 = I_1 \backslash \{ <\beta>:\beta < \alpha_{<>} \}$, \newline
$I_3 = \{ <\beta>:\beta < \alpha = \alpha_{<>} \}$. \newline
Note that by the inductive choice of the $J_k$'s:
\medskip
\roster
\item "{$\otimes$}"  if $\eta \in S_f \backslash S_f[I_1]$
\underbar{then} $\{ a^\ell
_{\eta,m}:\ell < n \text{ and } m < \omega \}$ 
has a finite intersection with $Y[I_1]$.
\endroster
\medskip

\noindent
(Use (3) noting that if $\eta \in S_f \backslash S_f[I_1]$ then
$\eta(\ell_k) > \alpha_{\eta \restriction \ell_k}$ for every
$k < k^*$ such that \newline
$\eta \restriction \ell_k \in J_k$). \newline
Note also that: $I_0 \in K,I_1 \in K^+,I_2 \in K^+,I_3 \in K^+$.  Also
$I_3 \le^d I_1$ and $I_3 \le^d I_0$ and $I_2 \le I_1 \le I_0$ (see 
Definition C(1)) and $G_{I_0} = G$. \newline

Note that $g_{I_0,I_1}$ is well defined
(see Definition J and Fact K). \newline
[Why?  We have to check $\otimes_{I_1,I_0}$ as
defined there, but $\otimes$ above says this].  Note also that $g_{I_1,I_2}$
is well defined (again we have to check $\otimes_{I_2,I_1}$ as defined in
Definition J, but for $\eta \in S_f(I_1) \backslash S_f(I_2)$ by their
definitions, $\eta(0) < \alpha_{<>}$ so easily $\dsize \bigcup_{\ell < n}
s^\ell_\eta$ is disjoint to the required set).
Look at the sequence $G = G_{I_0} \,\, {\underset g_{I_0,I_1}\to 
\longrightarrow}$
\,\, $G_{I_1} \,\, {\underset g_{I_1,I_2}\to \longrightarrow}$ \,\, $G_{I_2}$.

We know that $G_{I_2}$ is free (by Fact A as $I_2 \in K^+$),
$g_{I_1,I_2}$ is a
homomorphism from $G_{I_1}$ onto $G_{I_2}$
(see above, by Fact K, clause ($\alpha$) and $\otimes$ above)
hence $Ker(g_{I_1,I_2})$ is a direct summand of
$G_{I_1}$, so there is a projection $g^*$ of $G_{I_1}$ onto
$Ker(g_{I_1,I_2})$.
Also $h_{I_3,I_1},h_{I_3,I_0}$ are embeddings (by conclusion H)
as $I_3 <^d I_1,
I_3 <^d I_0$, (check or see above).  
Also $h_{I_3,I_1}(G_{I_3}) = Ker(g_{I_1,I_2})$
(compare
Fact G clause ($\gamma$) and Fact K clause ($\beta$)).  Hence
$h_{I_3,I_0} \circ h^{-1}_{I_3,I_1} \circ g^* \circ g_{I_0,I_1}$ is a
projection from $G = G_{<>} = G_{I_0}$ onto $\text{Rang}(h_{I_3,I_0})$ i.e.
essentially $G_{<>,\alpha}$.  This finishes the proof of the main fact,
hence the theorem 1.2. \hfill$\square_{1.2}$
\enddemo

[Question: here we can increase $\alpha_\eta$; can we make it exact?  (See
Appendix 3.6)].
\bigskip

\proclaim{1.3 Claim}  We can strengthen the conclusion of 1.2 to: for any
given $W \subseteq \lambda$ we can demand: there is a $\lambda$-free non-free
group $G$ with set of elements $\lambda$ such that

$$
\align
\{ \delta \in W:&G \restriction \delta \text{ is a subgroup of } G, \\
  &\text{ and is a free direct summand of } G\}
\endalign
$$
\medskip

\noindent
is a stationary subset of $\lambda$.
\endproclaim
\bigskip

\demo{Proof}  In the proof of 1.2;
\medskip
\roster
\item "{(A)}" for any $W_0 \subseteq \{ \alpha < \lambda:\langle \alpha
\rangle \in S\}$ stationary subset of $\lambda$, we can replace $S$ by
$\{ \eta:\eta \in S \text{ and } \ell g(\eta) > 0 \Rightarrow \eta(0) \in
W_0\}$
\item "{(B)}" assuming that the set of member of $G$ is $\lambda$ then 
\newline
$\{ \delta < \lambda:\delta \text{ is the set of elements of }
G_{<>,\delta}\}$ is a club of $\lambda$.
\endroster
\medskip

\noindent
Together with Main Fact L and Fact I, we are done. \hfill$\square_{1.3}$
\enddemo
\bigskip

\demo{1.4 Discussion}  We can rephrase the proof of 1.1 combinatorially; i.e.
explicitly write a set of generators $X$ such that $G = G_\alpha \oplus
\langle X \rangle_G$, do not think it is clearer.  To some extent this is done
in Fact A of the proof of 2.2.
\enddemo
\newpage

\head {\S2 The General Case: \, for a variety} \endhead
\bigskip

We note here that a parallel theorem holds for any suitable variety
considering two variants of $\lambda$-separable (see Definition 2.1(2) and
Definition 2.4).  We do
the general case in less details.
\bigskip

\definition{2.1 Definition}  1)  $T$ is a variety \underbar{if} $T$ is a
theory (in a vocabulary $\bold \tau$) all whose axioms are equations 
or just has the form $\forall x_1,\dotsc,
x_n\varphi$, $\varphi$ an atomic formula. Without loss of generality every
member of $\bold \tau$ (function symbol or predicate) appears in some axiom
of $T$. \newline
2)  A model $M$ of $T$ is called $\lambda$-separable \underbar{if} for every
$A \subseteq M,|A| < \lambda$ we can represent $M$ as a free product
$M_1 * M_2$ such that $A \subseteq M_1$ and $M_2$ is free. \newline
3)  $T$ has the $n$-th $h$-construction principle \underbar{if} we can find
$N$, $b_{\ell,m}$ \newline
(for $\ell < n,m < \omega$) and $N_{\bar m}$ (for $\bar m
\in {}^n \omega$) such that:
\medskip
\roster
\item "{(i)}"  $N$ a model of $T$ of cardinality $\le |T| + \aleph_0$
\item "{(ii)}"  $N$ is free, moreover, for each $\ell^* < n$ and
$m^* < \omega$ we can complete \newline
$\{ b_{\ell,n}:\ell < n,m < \omega$ and $[\ell = \ell^* \Rightarrow m < m^*]
\}$ to a free basis of $N$, call the set of additional elements $C_{\ell,m}$
{\roster
\itemitem{ $(iii)(\alpha)$ }  if $\bar m^i = \langle m^i_\ell:
\ell < n \rangle \in {}^n \omega$ (for $i = 1,2$) and
$\bar m^1 \le \bar m^2$ (i.e. \newline
$(\forall \ell < n)(m^1_\ell \le m^2_\ell)$ 
\underbar{then} $N_{\bar m^2} \subseteq N_{\bar m^1} \subseteq N$,
\itemitem{ $(\beta)$ }  $b_{\ell,m} \in N_{\langle m_k:k < n \rangle} 
\Leftrightarrow m \ge m_\ell$ and
\itemitem{ $(\gamma)$ }  $N$ is the free product $N_{\bar m} *
 \langle \{ b_{\ell,m}:\ell < n,m < m_\ell \} \rangle_N$.
\endroster}
\item "{(iv)}"  for no free model $F$ of $T$, is $N*F / \langle b_{\ell,m}:
\ell < n,m < \omega \rangle_N$ free (equivalently $N*F$ has a free basis
extending $\{ b_{\ell,m}:\ell < n,m < \omega \}$).
\endroster
\enddefinition
\bigskip

\remark{2.1A Remark}  On the $n$-th construction principle see Eklof Mekler
\cite{EM2} and then Mekler Shelah \cite{MkSh:366}.  The difference (between 
the $n$-th construction
principle and the $n$-th $h$-construction principle) is
clause (iii), it is not clarified here if it adds anything.  In all cases
the hope is that the analysis of \cite{Sh:161},\S3,\S4 exhausts the reasons
of the existence of the desired complicated object in $\lambda$, and
the crucial parameter of the system $S$ (see beginning of the proof of 1.2
or \S3) is $n = n(S)$.  So the hope is that for each $T$, the class of 
cardinals $\lambda$
where we have an example is, for some $\alpha^* \le \omega$

$$
\align
{\frak C}_{\alpha^*} = \biggl\{ \lambda:&\text{there are } n,S,
\langle \lambda(\eta,S):\eta \in S_i \rangle, \langle B_\eta:\eta \in
S_c \rangle \\
  &\langle s^\ell_\eta:\eta \in S_f,\ell < n \rangle \text{ as in 3.6, 3.7 
and } n < \alpha^* \biggr\}.
\endalign
$$
\medskip

\noindent
Usually we deal with varieties with countable vocabulary.
\endremark
\bigskip

\proclaim{2.2 Theorem}  Assume there is a $\lambda$-free not $\lambda^+$-free
abelian group exemplified by $n$, $S$, $\langle
s^\ell_\eta:\ell < n$ and $\eta \in S_f \rangle$ as in the proof of 
1.2 and the theory $T$ has the $n$-th
$h$-construction principle and $|T| < \lambda$.
\medskip

\underbar{Then} $T$ has a $\lambda$-separable model of cardinality $\lambda$
which is not free.
\endproclaim
\bigskip

\demo{2.2A Conclusion}  If there is a $\lambda$-free not $\lambda^+$-free
abelian group \underbar{then} for the variety of groups (not the abelian one)
there is a $\lambda$-free, $\lambda$-separable group $G$ of cardinality
$\lambda$ which is not free.  (I.e. $G$ is a non-free group of cardinality
$\lambda,G$ can be represented as $\dsize \bigcup_{\alpha < \lambda}
G_\alpha,G_\alpha$ increasing continuously of cardinality $< \lambda$, each
$G_\alpha$ free and $G$ is the free product (for the variety of groups) of
$G_{\alpha + 1}$ and some $H_{\alpha + 1}$ for each $\alpha < \lambda$).
\enddemo
\bigskip

\demo{Proof of 2.2A}  We should just check the condition of 2.1(3) which is
straight as in \cite{Sh:161}. \newline
[I.e. let $N$ be the group freely generated by \newline
$\{ b_{\ell,m}:\ell \in [1,m)
\text{ and } m < \omega\} \cup \{ y_m:m < \omega\}$, let:
\medskip
\roster
\item "{(a)}" $b_{0,0} =: y_0$
\item "{(b)}"  $b_{0,m+1}$ is the product $b_{1,m+1} \, b_{2,m+1} \, ... \,
b_{n-1,m+1} \, b_{0,m}(y_{0,m+1})^2$
\item "{(c)}"  $C_{\ell,m} = \{ y_k:k \in [m,\omega)\} \cup \{ b_{\ell,0}\}$
\item "{(d)}"  for $\bar m \in {}^n \omega$ clearly \newline
$\{ b_{\ell,m}:\ell \in [1,n) \text{ and } m < \omega\} \cup \{ b_{0,n}: n <
m_0\} \cup C_{0,m}$ is a free basis of $N$ and let $N_{\bar m}$ be the
subgroup of $N$ generated by \newline
$\{ b_{\ell,m}:\ell < n \text{ and } m \in
[m_\ell,\omega)\} \cup \{ y_m:m \in [m_0,\omega)\}$.
\endroster
\medskip

\noindent
Now check].
\enddemo
\bigskip

\demo{Proof of 2.2}  Let $\langle N,b_{\ell,m},N_{\bar m}:
\ell < n,m < \omega \text{ and } \bar m \in {}^n \omega \rangle$
exemplify the $n$-th $h$-construction principle.  We choose $n,S,...$ as in
the proof of 1.2.

Let $M$ be freely generated by $x[a]$ (for $a \in \dsize
\bigcup_{\eta \in S_c}
B_\eta)$ and $y_{\eta,c}$ (for $\eta \in S_f$ and $c \in N$) except that:
\medskip
\roster
\item "{(i)}"  $y_{\eta,c} = x[a]$ \underbar{if} $c = b_{\ell,m}$ and
$a = a^\ell_{\eta,m}$
\item "{(ii)}"  $\varphi(y_{\eta,c_1},\dotsc,y_{\eta,c_k})$
whenever $N \models ``\varphi(c_1,\dotsc,c_k)$" and $\varphi$ is a \newline
$\tau$-atomic formula.
\endroster
\enddemo
\bigskip

\subhead {Fact A} \endsubhead  For $\alpha < \lambda$ such that
$\langle \alpha \rangle
\notin S$ we can find $Y_0,Y_1,Y_2,S_0,S_1,S_2$ such that:
\medskip
\roster
\item "{(a)}"  $S_2 = S_f,Y_2 = \dsize \bigcup_{\eta \in S_c} B_\eta$
\item "{(b)}"  $S_0 = \{ \eta \in S_f:\eta(0) < \alpha \}$ and
$Y_0 = \cup \{ B_\eta:\eta \in S_c$ and $\eta(0) < \alpha \}$
\item "{(c)}"  $S_0 \subseteq S_1 \subseteq S_2$ and $Y_0 \subseteq Y_1
\subseteq Y_2$ and $Y_1$ is downward closed (remember $Y_2$ is a tree, see
3.6) so
$a^\ell_{\eta,m} \in Y_1 \and m_1 < m \Rightarrow a^\ell_{\eta,m_1} \in
Y_1$
\item "{(d)}"  for $\eta \in S_2 \backslash S_1$ the set $\{ a^\ell_{\eta,m}:
\ell < n,m < \omega \} \cap Y_1$ is finite
\item "{(e)}"  there is a list $\langle \eta_\zeta:\zeta < \zeta^* \rangle$
of $S_1 \backslash S_0$ without repetitions and $\langle \ell(\zeta):
\zeta < \zeta^* \rangle$ such that $0 \le \ell(\zeta) < n$ and
$\langle m(\zeta):\zeta < \zeta^* \rangle,m(\zeta) < \omega$ such that:
{\roster
\itemitem{ $(\alpha)$ }  $\{ a^{\ell(\zeta)}_{\eta_\zeta,m}:m \in
[m(\zeta),\omega) \}$ is disjoint to \newline
$Y_0 \cup \{ a^\ell_{\eta_\varepsilon,
m}:\ell < n,\varepsilon < \zeta,m < \omega \}$
\itemitem{ $(\beta)$ }  $\{ a^{\ell(\zeta)}_{\eta_\zeta,m}:m < \omega \}
\subseteq Y_1$.
\endroster}
\endroster
\bigskip

\demo{Proof}  Included in the proof of Theorem 1.2.
\enddemo
\bigskip

\definition{Remark B}  We can add
\medskip
\roster
\item "{(f)}"  $S_1$ is $S_f[I_1]$ from the proof of Theorem 1.2 so for some
function $k$ from $S_1 \backslash S_0$ to $n = \{ 0,\dotsc,n-1\}$ we have
\newline
$Y_1 = Y_0 \cup 
\{ a^\ell_{\eta_\zeta,m}:\eta \in S_1,m < \omega \text{ and }
\ell \in [k(\eta),n)\}$.
\endroster
\enddefinition
\bigskip

\subhead {Fact C} \endsubhead  Under the conclusion of Fact A, letting
\newline
$M_0 =: \langle \{ x[a]:a \in Y_0 \} \cup \{ y_{\eta,c}:\eta \in S_0,
c \in N \} \rangle_M$ we have: $M_0$ is free and for \newline
some $M_2$, $M = M_0*M_2$.
\bigskip

\demo{Proof}  Clearly $M_1$ is free (for $T$) as in the proof of Fact A in
the proof of 1.2.  The new point is to find $M_2$. \newline
For each $\ell < n,m < \omega$, let $C_{\ell,m} \subseteq N$ be such that
\newline
$C_{\ell,m} \cup \{ b_{\ell_1,m_1}:\ell_1 \ne \ell,m_1 < \omega$ \underbar
{or} $\ell_1 = \ell,m_1 < m \}$ is a free basis of $N$ with no
repetitions.
\medskip

We let $M_2$ be the submodel of $M$ generated by:
\medskip
\roster
\item "{(A)}"  $y_{\eta_\zeta,c}$ \underbar{if} $\zeta < \zeta^*$ and
$c \in C_{\ell(\zeta),m(\zeta)}$
\item "{(B)}"  $x[a^\ell_{\eta_\zeta,m}]$ \underbar{if} $\zeta < \zeta^*,
a^\ell_{\eta_\zeta,m} \in Y_1 \backslash Y_0$ and for no $\varepsilon <
\zeta^*$ do we have \newline
$a^\ell_{\eta_\zeta,m} \in \{a^{\ell(k)}_{\eta_\varepsilon,k}:
k \in [m(\varepsilon),\omega) \}$
\item "{(C)}"  $x[a]$ \underbar{if} $a \in Y_2 \backslash Y_1$
\item "{(D)}"  $y_{\eta,c}$ \underbar{if} $\eta \in S_2 \backslash S_1,
c \in N_{\langle m_\ell(\eta):\ell < n \rangle}$ where
$m_\ell(\eta) = \text{ min} \{m:a^\ell_{\eta,m} \notin Y_1 \}$.
\endroster
\medskip

\noindent
Now
\medskip
\roster
\item "{$(*)_1$}"  $M = \langle M_0,M_2 \rangle$.
\endroster
\medskip

First we prove by induction on $\zeta < \zeta^*$ that $\{x[a^\ell
_{\eta_\zeta,m}]:\ell < n \text{ and } m < \omega\} \subseteq \langle
M_0,M_2 \rangle$ and $\{ y_{\eta_\zeta,c}:c \in N\} \subseteq \langle M_0,
M_2 \rangle$.  Arriving to $\zeta$ we split the proof to cases.
\bigskip

\noindent
\underbar{Case 1}:  $a^\ell_{\eta_\zeta,m} \in Y_0$. \newline
Then $x[a^\ell_{\eta_\zeta,m}] \in Y_0 \subseteq M_0 \subseteq \langle
M_0,M_2 \rangle$.
\bigskip

\noindent
\underbar{Case 2}:  $a^\ell_{\eta_\zeta,m} \in Y_1 \backslash Y_0$ and for
some $\varepsilon < \zeta,a^\ell_{\eta_\zeta,m} \in \{a^{\ell(\varepsilon)}
_{\eta_\varepsilon,k}:k \in [m(\varepsilon),\omega)\}$. \newline
We use the induction hypothesis on $\varepsilon$.
\bigskip

\underbar{Case 3}:  $a^\ell_{\eta_\zeta,m} \in Y_1 \backslash Y_0$ and 
$\ell \ne \ell(\zeta) \vee [\ell = \ell(\zeta) \and m < m(\zeta)]$ 
\newline
and for no $\varepsilon < \zeta$, do we have 
$a^\ell_{\eta_\zeta,m} \in \{a^{\ell(\varepsilon)}
_{\eta_\varepsilon,k}:k \in [m(\varepsilon),\omega)\}$. \newline

Now 
$\varepsilon < \zeta^*$ implies $a^\ell_\eta \notin \{ a^{\ell(\varepsilon)}
_{\eta_\varepsilon,k}:k \in [m(\varepsilon),w)\}$. \newline
[Why?  If $\varepsilon < \zeta$ this is assumed in the case, if $\varepsilon
= \zeta$ this is follows by $\ell \ne \ell(\zeta)$, and if $\varepsilon \in
(\zeta,\zeta^*)$ this follows by clause $(e)(\alpha)$ (with $\varepsilon$'s
here standing for $\zeta,\varepsilon$ there).  Hence the assumption of
clause $(B)$ holds.]

By clause $(B)$, $x[a^\ell_{\eta_\zeta,m}] \in M_2 \subseteq \langle M_0,
M_2 \rangle$.
\bigskip

\underbar{Case 4}:  $a^\ell_{\eta_\zeta,m} \in Y_2 \backslash Y_1$. \newline
By clause $(C)$, $x[a^\ell_{\eta_\zeta,m}] \in M_2 \subseteq \langle M_0,
M_2 \rangle$.
\bigskip

\underbar{Case 5}:  No previous cases. \newline

By the earlier cases $\ell = \ell(\zeta)$ and \newline
$\{ x[a^{\ell^*_1}_{\eta_\zeta,m_1}]:\ell^*_1 < n,m^*_1 < \omega \text{ and }
[\ell^*_1 \ne \ell(\zeta) \Rightarrow m^*_1 < m(\zeta)]\} \subseteq \langle 
M_0,M_2 \rangle$.
\medskip

Let $N' =: \{ c \in N:y_{\eta_\zeta,c} \in \langle M_0,M_2 \rangle\}$, so by
the previous sentence \newline
$\{b_{\ell_1,m_1}:\ell_1 < n,m_1 < \omega \text{ and }
\ell_1 = \ell(\zeta) \Rightarrow m_1 < m(\zeta)\} \subseteq N'$, and by
clause $(A)$ \newline
also $C_{\ell(\zeta),m(\zeta)} \subseteq N'$ hence (see
clause $(ii)$ in Definition 2.1) clearly $N' = N$, \newline
so $x[a^{\ell_1}_{\eta_\zeta,m}] \in \langle M_0,M_2 \rangle$ and 
$y_{\eta_\zeta,c} \in \langle M_0,M_2 \rangle$. \newline

We have proved $\{ x[a]:a \in Y_1 \backslash Y_0\} \subseteq 
\{x[a^\ell_{\eta_\zeta,
m}]:\ell < n,m < \omega,\zeta < \zeta^*\} \subseteq \langle M_0,M_2 \rangle$.
As $\{x[a]:a \in Y_0\} \subseteq M_0 \subseteq \langle M_0,M_2 \rangle$ and
by clause $(C)$ we have $\{x[a]:a \in Y_2 \backslash Y_1 \} \subseteq
\langle M_0,M_2 \rangle$ we conclude $\{x[a]:a \in Y_2 \} \subseteq
\langle M_0,M_2 \rangle$.

Also we have proved $\{y_{\eta_\zeta,c}:c \in N,\zeta < \zeta^*\} \subseteq
\langle M_0,M_2 \rangle$ (this was done during the proof of case 5) so 
$\{y_{\eta,c}:\eta \in S_1 \backslash S_0
\text{ and } c \in N\} \subseteq \langle M_0,M_2 \rangle$. \newline
Also for $\eta \in S_2 \backslash S_1$, letting $N^\eta = \{c \in N:
y_{\eta,c} \in \langle M_0,M_2 \rangle\}$, \newline
$m_\ell = \text{ min}\{m:a^\ell
_{\eta,m} \notin Y_1 \}$ we have: by clause $(D)$, $N_{\langle m_\ell:\ell
< n \rangle} \subseteq N^\eta$, and \newline
$\{ a^\ell_{\eta,m}:\ell < n,m < \omega
\} \subseteq Y_2 \text{ so } b_{\ell,m} \in N^\eta \text{ hence }
N^\eta = N$ so \newline
$\{ y_{\eta,c}:\eta \in S_2 \backslash S_1,c \in N\} \subseteq
\langle M_0,M_2 \rangle$.  Lastly if $\eta \in S_0$ we have \newline
$\{ y_{\eta,c}:c \in N\} \subseteq M_0 \subseteq \langle M_0,M_2 \rangle$.
Together
$\{ y_{\eta,c}:\eta \in S_2 \text{ and } c \in N\} \subseteq \langle 
M_0,M_2 \rangle$; and also we note above $\{x[a]:a \in Y_2\} \subseteq
\langle M_0,M_2 \rangle$; we can conclude $M = \langle M_0,M_2 \rangle$, i.e.
$(*)_1$.
\bigskip

\noindent
So to finish the proof we need
\medskip
\roster
\item "{$(*)_2$}"  $M = M_0 * M_2$
\endroster
\medskip

\noindent
(i.e. they generate $M$ freely). \newline 
Look at the equations in the definition of $M$ and together with the proof
of $(*)_1$ rewrite them in terms of
the generators of $M_0$ and of $M_2$.  The equations either trivialized or
speak on generators of $M_0$ or speak on generators of $M_2$.  [more?]
\hfill$\square_{2.2}$
\enddemo
\bigskip

\noindent
Note that as the variety of abelian groups is very nice, e.g. a subgroup
of a free abelian group is free, distinct definitions for general varieties
become identified for it; so Theorem 1.2 has various generalizations and
Theorem 2.2 is not the unique one.  Another generalization is presented below.
\proclaim{2.3 Theorem}  Assume $\lambda$ is as in 1.2 with $n,S,\langle s^\ell
_\eta:\ell < n,\eta \in S_f \rangle$ such that $T$ has the $n$-th construction
principle (i.e. in Definition 2.1 we omit clause (iii), but demanding each
$C_{\ell,n}$ is infinite; this holds without loss of generality by clause
(iv) of Definition 2.1).  \underbar{Then} there is a model $M$ of $T$, not
free of cardinality $\lambda$, but is $\lambda$-proj-separable, where:
\endproclaim
\bigskip

\definition{2.4 Definition}  For a variety $T$ and a model $M$ of $T$ and
cardinality $\lambda$ we say $M$ is $\lambda$-proj-separable, if for
every $A \subseteq M,|A| < \lambda$
there is a free $M' \subseteq M$ including $A$ and a projection $h$ from $M$
onto $M'$.
\enddefinition
\bigskip

\demo{Proof of 2.3}  We define $M,x[a],y_{\eta,c}$ as in the proof of 2.2.
For every $\ell(*) < n$ and $m(*) < \omega,m(*) > 0$ there is a homomorphism 
$g_{\ell(*),m(*)}$ from $N$ onto 
$\langle b_{\ell,n}:\ell < n,m < \omega \text{ and } [\ell = \ell(*) 
\Rightarrow m < m(*)] \rangle_N$ which is the identity on
$\langle b_{\ell,m}:\ell < n,m < \omega$ and $[\ell = \ell(*) \Rightarrow
m < m(*)] \rangle_N$ (maps the members of $C_{\ell(*),m(*)}$ onto
$\{ b_{\ell(*),0} \}$.)
\newline
Let $\Gamma$ be the set of equations which we make the generators satisfy.
We choose $Y_0,Y_1,Y_2,S_0,S_1,S_2$ as in Fact A from the proof of 2.2 and
without loss of generality $\zeta < \zeta^* \Rightarrow m(\zeta) > 0$.  Let
$\{ \eta_\zeta:\zeta \in [\zeta^*,\zeta^{**})\}$ list $S_2 \backslash S_1$.

For each $\eta_\zeta \in S_2 \backslash S_1$ we can choose

$$
\ell(\zeta) = 
\ell_{k^*-1},m(\zeta) = \text{ Min}\{m:0 < m < \omega \text{ and }
a^{\ell(\zeta)}_{\eta,m} \notin Y_1\}.
$$
\medskip

\noindent
Let $M_1$ be the model of $T$ generated by \newline
$\Xi_1 = \{x[a^\ell_{\eta,m}]:\ell < n,\eta \in S_1,m < \omega\} \cup
\{ y_{\eta,c}:c \in N,\eta \in S_2\}$ freely except

$$
\align
\Gamma_1 = &\biggl\{ y_{\eta,c} = x[a]:c = b_{\ell,m},a = a^\ell_{\eta,m}
\text{ and } \eta \in S_1 \biggr\} \cup \\
  &\biggl\{ \varphi(y_{\eta,c_1},\dotsc,y_{\eta,c_k}):N \models \varphi
(c_1,\dotsc,c_k),\varphi \text{ a } T \text{-atomic formula}\biggr\}.
\endalign
$$
\medskip

\noindent
Let $M^-_2$ be the model of $T$ generated by (note: $I_1,J_k$ are from the
proof of 1.2) 

$$
\align
\Xi^-_2 =: \biggl\{ x[a]:&a \in Y_2 \text{ but if } a \in B_{\eta \char 94
\langle \lambda(\eta,S) \rangle}, \\
  &\ell g(\eta) = \ell_{k^*-1} \text{ and } \eta \in J_{k^*-1} \subseteq 
J[I_1] \text{ then } a \text{ is in the first level} \\
  &\text{(i.e. like } a^{\ell_{k^*}}_{\eta,0}) \text{ or } a \in
B_{\eta \char 94 \langle \alpha_\eta \rangle} \\
  &(\alpha_\eta \text{ from the choice of } I_1) \biggr\} 
\endalign
$$   

\noindent
freely except the equations

$$
\align
\Gamma^-_2 = \Gamma_1 = &\biggl\{ y_{\eta,c} = x[a]:
c = b_{\ell,m},a = a^\ell_{\eta,m},\eta \in S_1,\ell < n, m < \omega 
\text{ and } x[a] \in \Xi^-_2 \biggr\} \cup \\
  &\biggl\{ \varphi(y_{\eta,c_1},\dotsc,y_{\eta,c_k}):N \models \varphi
[c_1,\dotsc,c_k],\varphi \text{ a } T \text{-atomic formula},
\eta \in S_1 \biggr\}.
\endalign
$$
\medskip

\noindent
(Note that if $\eta \in J_{k-1}$ and $\eta \triangleleft \nu \in S_f$ then 
cf$(\nu(\ell_{k^*-1})) = \aleph_0$). \newline
Clearly $M_0 \subseteq M_1 \subseteq M^-_2 \subseteq M$. \newline
We define a homomorphism $h_2$ from $M$ into $M^-_2:h_2 \restriction M^-_2$
is the identity, and for $\eta = \eta_\zeta \in S_2 \backslash S_1$ and
$c \in N$ we let:

$$
h_2(y_{\eta,c}) = y_{\eta,g_{\ell(\zeta),m(\zeta)}(c)}.
$$
\medskip

\noindent
Note: $h_2(x[a^\ell_{\eta_\zeta,m}]) = x[a^\ell_{\eta_\zeta,m}]$ when
$\ell \ne \ell(\zeta) \vee m < m(\zeta)$ by the tree structure of
$\dsize \bigcup_{\eta \in S_c} B_\eta$, the cases of the definition of
$h_2$ are compatible and the equations are preserved.  So $h_2$ is a
homomorphism and even a projection from $M$ onto $M^-_2$. \newline

Trivially, we can find a projection $h_1$ from $M^-_2$ onto $M_1$.

Next note that $M_1$ is a free extension of $M_0$ (a free basis is
\newline
$\{y_{\eta_\zeta,c}:c \in C_{\ell(\zeta),m(\zeta)} \text{ and } \zeta <
\zeta^*\} \cup \{x[a]:a \in Y_1 \backslash Y_0 \text{ and for no }
\zeta < \zeta^* \text{ is}$ \newline
$a \in \{ a^{\ell(\varepsilon)}_{\eta_\zeta,m}:
m \in [m(\varepsilon),\omega)\}\}$.
\medskip

So we can find a projection $h_0$ from $M_1$ onto $M_0$.  So $h_0 \circ
h_1 \circ h_0$ is a projection as required. \hfill$\square_{2.3}$
\enddemo
\bigskip

\proclaim{2.5 Claim}  Theorems 2.2, 2.3 can be strengthened as in 1.3.
\endproclaim
\bigskip

\demo{2.6 Discussion}  Implicit in the proof of 2.3 is an alternative
criterion sufficient for the conclusion of 2.2.
\enddemo
\newpage

\head {\S3 Appendix: characterizing the existence in $\lambda$ of an almost
free abelian group} \endhead
\bigskip

To make the main theorem 1.2 more easily read we present part of
\cite{Sh:161}, more exactly a variant to \cite{Sh:161},3.6,3.7,p.212.
\newline
Numbers are as in \cite{Sh:161}.
\bigskip

\definition{3.1 Definition} \newline
(1) For a regular uncountable cardinal $\lambda(> \aleph_0)$ we call $S$ a
$\lambda$-set if:
\medskip
\roster
\item "{(a)}"  $S$ is a set of strictly decreasing sequences of ordinals
$< \lambda$.
\item "{(b)}"  $S$ is closed under initial segments and is nonempty.
\item "{(c)}"  For $\eta \in S$ if we let $W(\eta,S) =: \{i:\eta \char 94
<i> \in S \}$ and \newline
$\lambda(\eta,S) =: \text{ Sup }W(\eta,S)$ \underbar{then} 
whenever $W(\eta,S)$ is not empty,
$\lambda(\eta,S)$ is a regular uncountable cardinal and
$W(\eta,S)$ is a stationary subset of $\lambda(\eta,S)$.  
Also $\lambda(<>,S) = \lambda$ (and by clause (a) we know
$\lambda(\eta \char 94 < \alpha >,S) \le | \alpha |)$.
\endroster
\medskip

\noindent
(2) For a $\lambda$-set $S$, let $S_f(=$ set of final elements of $S$) be
$\{\eta \in S:(\forall i) \eta \char 94 <i> \notin S \}$ and
$S_i(=$ set of initial elements of $S$) be $S \backslash S_f$ so $(S_f = \{ \eta \in S:
\lambda(\eta,S) = 0 \}$).

We sometimes allow $\lambda = 0$.  Then the only $\lambda$-set is $\{ <> \}$.
\newline
(3) For $\lambda$-sets $S^1,S^2$ we say $S^1 \le S^2$ \, ($S^1$ a sub-$\lambda$-set of $S^2$) if $S^1 \subseteq S^2$ and $\lambda(\eta,S^1) = \lambda(\eta,
S^2$) for every $\eta \in S^1$ (so $S^1_i = S^1 \cap S^2_i$).  Clearly
 $\le$ is transitive.
\enddefinition
\bigskip

\noindent
Notation: In this section $S$ will be used to denote $\lambda$-sets.
\remark{3.1A Remark}  Many of the properties below holds also if we waive the
``decreasing" demand in clause $(a)$ but not all, and for what we want to
analyze we can get such $S$, so we have concentrated on this family of sets.
\endremark
\bigskip

\proclaim{3.2 Claim} \newline
(1) If $S$ is a $\lambda$-set, $\lambda(\eta,S) > \kappa$ for every
$\eta \in S_i$ (holds always for $\kappa = \aleph_0$) and $G$ is a function
from $S_f$ to $\kappa$ \underbar{then} for some $S^1 \le S$ we have:
$G$ is constant on $S^1_f$. \newline
(2)  If $S$ is a $\lambda$-set and $\eta \in S_i$, \underbar{then}
$S^{[\eta]} = \{ \nu:
\eta \char 94 \nu \in S \}$ is a $\lambda(\eta,S)$-set, and \newline
$\lambda(\nu,S^{[\eta]}) = \lambda(\eta \char 94 \nu,S)$ for every
$\nu \in S^{[\eta]}$. \newline
(3)  If $S$ is a $\lambda$-set, $\kappa$ a regular cardinal, $(\forall \eta
\in S) \, (\lambda(\eta,S) \ne \kappa)$ and $G$ is a function from $S$ to
$\kappa$ \underbar{then} for some $S^1 \le S$ and $\gamma < \kappa$ for every
$\eta \in S^1$ we have $G(\eta) < \gamma$. \newline
(4) If $\lambda > \aleph_0$ is regular, $W \subseteq \lambda$ stationary, for
$\delta \in W,\, S^\delta$ is a $\lambda_\delta$-set
for some $\lambda_\delta \le \delta$ or $S^\delta = \{ <> \}$ \underbar{then}
$S =: \{ <> \} \cup \{
 \langle \delta \rangle \char 94
\eta:\eta \in S^\delta,\delta \in W \}$ is a $\lambda$-set and 
$\lambda(<\delta>
\char 94 \eta,S) = 
\lambda(\eta,S^\delta)$ for $\delta \in W,\eta \in S^\delta$ and
$S_i = \{ \langle \rangle \} \cup \dsize \bigcup_{\delta \in W} S^\delta_i$.
\newline
(5)  If $S$ is a $\lambda$-set, $F$ a function with domain $S \backslash
\{ <> \}$,
$F(\eta \char 94 \langle \alpha \rangle)
 < 1 + \alpha$ \underbar{then} $F$ is essentially constant for
some $S^1 \le S$ which means $F \restriction \{ \eta \in S^1:\ell g(\eta) = 
m\}$ is constant for each $m$. \newline
(6)  For any $\lambda$-set $S$ there is a $\lambda$-set $S^1 \le S$ such that:
\medskip
\roster
\item "{(a)}"  all $\eta \in S_f$ has the same length $n$
\item "{(b)}"  for each $\ell < n$ either
{\roster
\itemitem{ (i) }  every $\eta(\ell) \, (\eta \in S_f)$ is an inaccessible
cardinal (not necessarily strong limit), \underbar{or}
\itemitem{ (ii) }  every $\eta(\ell) \, (\eta \in S_f)$ is a singular limit
ordinal,
\endroster}
\item "{(c)}"  for each $\ell < n$, either
{\roster
\itemitem{ (i) }  $\lambda(\eta \restriction (\ell + 1),S) = \eta(\ell)$ for
every $\eta \in S_f$ \underbar{or}
\itemitem{ (ii) }  $\lambda(\eta \restriction (\ell + 1),S) = \lambda^{\ell+1}
_S$ for every $\eta \in S_f$ (for a fixed $\lambda^{\ell + 1}_S$).
\endroster}
\item "{(d)}"  The truth value of
``$\text{cf}(\eta(\ell)) = \lambda(\eta \restriction m,S)$" 
is the same for all $\eta \in S_f$ (for constant $\ell,m < n$).
\endroster
\endproclaim
\bigskip

\demo{Proof}  Straightforward, e.g. \newline
(5)  In first glance we get only: if $\rho \in S_i$ then
$F \restriction \{ \rho \char 94 \langle \alpha \rangle:
\alpha \in W(\rho,S)\}$ is constant (by Fodor's lemma and the demand 
``$W(\rho,S)$ is a stationary subset of $\lambda(\rho,S)$".  However, as
every $\eta \in S$ is (strictly) decreasing sequence of ordinals we can
iterate this (simpler if we first apply part (6) clause $(a)$). 
\hfill$\square_{3.2}$
\enddemo
\bigskip

\proclaim{3.3 Claim}  Suppose $P$ is a family of sets which exemplify the
failure of $PT(\lambda,\kappa^+)$ (where $\lambda > \kappa$) i.e.
$a \in P = |a| \le \kappa,P$ has no transversal ($=$ one to one choice
function) but every $P' \subseteq P$ of cardinality $< \lambda$ has a
transversal.  \underbar{Then} there is a
$\lambda$-set $S$ and function $F$ with domain $S_f$ such that:
\medskip
\roster
\item "{(a)}"  For each $\eta \in S_f,F(\eta)$ is a subfamily of $P$ of power
$\le \kappa$.
\item "{(b)}"  For $\eta \in S_i$ we have $\lambda(\eta,S) > \kappa$.
\item "{(c)}"  For $\eta \in {}^{\omega >}(\lambda + 1)$, let $F^0(\eta)
 =: \cup
\{ F(\tau):\tau <_{\ell x} \eta$ and $\tau \in S_f \}$, 
where $<_{\ell x}$ is the
lexicographic order, $F^1(\eta) =: \cup \{ F(\tau):\eta \trianglelefteq
 \tau \in S_f \}$
and \newline
$F^2(\eta) =: \cup \{ A:A \in F^0(\eta \char 94 \langle \lambda(\eta,S)
\rangle)\} \backslash \cup \{ A:A \in F^0(\eta)\}$.
\endroster
\medskip

\noindent
Note that for $\eta \in S$ we have $F^2(\eta \char 94 \langle \lambda(\eta,S)
\rangle) = F^0(\eta) \cup F^1(\eta)$. \newline
\medskip
\roster
\item "{(d)}" ${{}}$
{\roster
\itemitem{ $(\alpha)$ }  For $\eta \in S_f,F^1(\eta)/F^0(\eta)$ 
is not free, (that
is $F^1(\eta)$ has no one to one choice function with range disjoint to
$\cup \{ A:A \in F^0(\eta) \}$).
\itemitem{ $(\beta)$ }  For $\eta \in S_i,F^1(\eta)/F^0(\eta)$ 
is $\lambda(\eta,S)$-free not free
and $|F^1(\eta)| = \lambda(\eta,S)$ (this follows as
$| \{ \tau:\eta \trianglelefteq \tau \in S \}|) = \lambda(\eta,S)$).
\endroster}
\item "{(e)}"  If $\eta \char 94 <\alpha> \, \in S$ \underbar{then}
$\alpha$ is a limit
ordinal, $\text{cf}(\alpha) \le \lambda(\eta \char 94 <\alpha>,S) + \kappa \le
|\alpha|$ and if $\beta < \lambda(\eta,S)$ is an inaccessible cardinal
$(> \aleph_0)$ then $\beta \cap W(\eta,S)$ is not a stationary subset of
$\beta$.
\item "{(f)}"  If $\eta \char 94 <\alpha> \triangleleft \,
 \nu \in S_f$ and $\text{cf}(\alpha) > \kappa$
then for some natural number $k$ we have $\eta \char 94 <\alpha> 
\trianglelefteq \nu \restriction k$ and
$\lambda(\nu \restriction k,S) = cf(\alpha)$ 
(so if $\alpha$ is an inaccessible cardinal then $k = \ell g(\eta)$).
\endroster
\endproclaim
\bigskip

\demo{Proof} See \cite{Sh:161}.
\enddemo
\bigskip

\remark{Remark}  Note clause (f), it is crucial; without it we won't be
able to prove the desired conclusion.
\endremark
\bigskip

\definition{3.4 Definition} \newline
(1)  A $\lambda$-system is $\langle B_\eta:\eta \in S_c \rangle$ where:
\medskip
\roster
\item "{(a)}"  $S$ is a $\lambda$-set, and we let $S_c =: \{ \eta \char 94
\langle i \rangle: \eta \in S_i$ and $i < \lambda(\eta,S) \}$
\item "{(b)}"  $B_{\eta \char 94 \langle i \rangle} \subseteq 
B_{\eta \char 94 \langle j \rangle}$ when $\eta \in S_i$ and $i < j < 
\lambda(\eta,S)$
\item "{(c)}"  If $\delta$ is a limit ordinal $< \lambda(\eta,S)$ 
\underbar{then}
$B_{\eta \char 94 \langle \delta \rangle} = \cup \{ B_{\eta \char 94
\langle i \rangle}:i < \delta \}$
\item "{(d)}"  $|B_{\eta \char 94 \langle i \rangle}| < \lambda(\eta,S)$ for
$i < \lambda(\eta,S)$.
\endroster
\medskip

\noindent
(2)  The $\lambda$-system $\langle B_\eta:\eta \in S_c \rangle$ is called
disjoint if the sets $\{ B_{\eta \char 94 \langle \lambda(\eta,S) \rangle}:
\eta \in S_i \}$ (see (3) below) are pairwise disjoint. \newline
(3)  We let $S_m =: S \backslash \{ <> \}$, 
$B_{\eta \char 94 \langle \lambda(\eta,S)
\rangle} =: B^*_\eta =: \cup \{ B_{\eta \char 94 \langle i \rangle}:
i < \lambda(\eta,S) \}$ for $\eta \in S_i$.
\enddefinition
\bigskip

\proclaim{3.5 Claim}  Suppose $\lambda$ is a regular uncountable cardinal,
$\langle B_\eta:\eta \in S_c \rangle$ a $\lambda$-system, and for $\eta \in
S_f, s_\eta \subseteq \dsize \bigcup_{\ell < \ell(\eta)} B_{\eta \restriction
(\ell + 1)}$.  \underbar{Then} $ \{s_\eta:\eta \in S_f \}$ has no transversal.
\endproclaim
\bigskip

\demo{Proof}  Straightforward (or see \cite{Sh:161}).
\enddemo
\bigskip

\proclaim{3.6 Claim}  Suppose $PT(\lambda,\kappa^+)$ fails (see 3.3).
\footnote{we are interested mainly in the case $\kappa = \aleph_0$} 
\underbar{Then} 
there is a disjoint $\lambda$-system \newline
$\langle B_\eta:\eta \in S_c \rangle$ and sets
$s^\ell_\eta$ (for $\eta \in S_f$ and $\ell < \ell g(\eta))$, and
$C_\delta$ (for $\delta < \lambda$ a limit ordinal) and 
$\varepsilon_{\eta,\ell}$ (for
$\eta \in S$ and $\ell < \ell g(\eta)$) such that:
\medskip
\roster
\item "{(a)}"  $S$ satisfies the conclusion
of Claims 3.2(6),3.3(e) and 3.3(f), in particular
$\eta \in S_f \Rightarrow \ell g
(\eta) = n$.
\item "{(b)}"  $s^\ell_\eta \subseteq B_{\eta \restriction(\ell + 1)},0 <
|s^\ell_\eta| \le \kappa$.
\item "{(c)}" For every $I \subseteq S_f$: \underbar{if} $|I| < \lambda$ 
\underbar{then}
$\{ \dsize \bigcup_\ell s^\ell_\eta:\eta \in I \}$ has a transversal
(as as indexed set).
Moreover, for every $\rho \in S_i$ if $I \subseteq \{ \nu:\rho
\trianglelefteq \nu \in
S_f \}$ and $|I| < \lambda(\rho,S)$ \underbar{then} the family \newline
$\{ \dsize \bigcup_{\ell
 \ge \ell g(\rho)} s^\ell_\eta:\eta \in I \}$ has a transversal.
\item "{(d)}"  If $s^\ell_\eta \cap s^m_\nu \ne \emptyset$ \underbar{then}
{\roster
\itemitem{ $(\alpha)$ }  $\ell = m$ and the sequences $\eta,\nu$ are 
different only at the $\ell$-th place i.e.
$\rho =: \eta \restriction \ell = \nu \restriction \ell$ and 
$\eta \restriction [\ell + 1,n) = \nu \restriction 
[\ell + 1,n)$ \underbar{and} 
\itemitem{ $(\beta)$ }  $\lambda(\eta \restriction i,S) = \lambda(\nu
\restriction i,S)$ when $\ell + 1 < i < n$ and
\itemitem{ $(\gamma)$ }  \underbar{either} $\lambda(\eta \restriction
(\ell + 1),S) = \eta(\ell)$ and $\lambda(\nu \restriction (\ell + 1),S) =
\nu(\ell)$ are both inaccessible cardinals \underbar{or} 
$\lambda(\eta \restriction
(\ell + 1),S) = \lambda(\nu \restriction (\ell + 1),S)$.
\endroster}
\item "{(e)}"  For $\eta \char 94 < \delta > \in S$ we have
{\roster
\itemitem{ $(\alpha)$ }  $C_\delta$ is a closed 
unbounded subset of $\delta$, $C_\delta
= \{ \zeta(\delta,i):i < \text{cf}(\delta) \}$, $\zeta(\delta,i)$ increasing
continuously with $i$
\itemitem{ $(\beta)$ }  In addition \underbar{if} $\nu = \eta 
\restriction \ell,\nu \in S_i,\eta \in S_i,
\lambda(\eta,S) = \text{cf}[\eta(\ell)] > \aleph_0$ \underbar{then}
$\varepsilon_{\eta,\ell}$ is a
strictly increasing function from $\lambda(\nu,S)$ to $\lambda(\nu,S)$
\itemitem{ $(\gamma)$ }  in clause $(\beta)$ if $\delta =: \eta(\ell)$ is an
inaccessible cardinal (hence necessarily $\ell g(\eta) = \ell + 1$) then
\newline
$\emptyset = W(\nu,S) \cap \{ \zeta(\delta,i):i$ belong to the range of
$\varepsilon_{\eta,\ell}\}$
\endroster}
\item "{(f)}" ${{}}$
{\roster
\itemitem{ $(\alpha)$ }  If $\ell < m < n,\eta \in S_f$, $\text{cf}[\eta
(\ell)] = \lambda(\eta \restriction m,S) > \kappa$ then \newline
$s^\ell_\eta \subseteq B_{(\eta \restriction \ell) \char 94 \langle \zeta + 1
\rangle} \backslash
B_{(\eta \restriction \ell) \char 94 \langle \zeta \rangle}$
where $\zeta = \zeta(\eta(\ell),\varepsilon_{\eta,\ell}(\eta(m))+2)$; i.e.
$\zeta$ is the $(\varepsilon_\eta(\eta(m)) + 2$)-th member of 
$C_{\eta(\ell)}$.  Moreover if
$s^\ell_\eta \cap s^\ell_\nu \ne \emptyset$, $\eta \ne \nu$ then
$\zeta(\eta(\ell),\eta(m)) = \zeta(\nu(\ell),\nu(m))$.
\itemitem{ $(\beta)$ }  If $\ell < m < n = \ell g(\eta),
\eta \in S_f,cf[\eta(\ell)] = \lambda(\eta \restriction m,S) \le \kappa$ 
\underbar{then} $s^\ell_\eta \subseteq
B_{\eta \restriction (\ell + 1)} \backslash
B_{\eta \restriction \ell \char 94 \langle \zeta
\rangle}$ where $\zeta = \zeta(\eta(\ell),\eta(m))$;
i.e. $\zeta$ is the \newline
$(\eta(m)+1)$-th member of $C_{\eta(\ell + 1)}$ and
$\xi < \eta(\ell) \Rightarrow |s^\ell_\eta \backslash
B_{(\eta \restriction \ell) \char 94 \langle \xi \rangle}| = \kappa$.
Moreover if $s^\ell_\eta \cap s^\ell_\nu \ne \emptyset,\eta \ne \nu$ then
\newline
$\zeta(\eta(\ell),\eta(m)) = \zeta(\nu(\ell),\nu(m)$).
\endroster}
\item "{(g)}"  If $\ell < \ell g(\eta), \eta \in S_f,cf[\eta(\ell)] \le
\kappa$ \underbar{then} for no $\zeta < \eta(\ell)$ is $s^\ell_\eta \subseteq
B_{\eta \restriction \ell \char 94 \langle \zeta \rangle}$.
\item "{(h)}"  For some well ordering $<^*_\eta$ of $B^*_\eta$ \,
$(\eta \in S_i)$ if $\eta \char 94 \langle i \rangle \trianglelefteq
\nu \in S_f$, \underbar{then} \newline
$[cf(i) \ge \kappa \Rightarrow s^{\ell g(\eta)}_\nu$ has order type $\kappa$]
and $cf(i) < \kappa \Rightarrow s^{\ell(\eta)}_\nu$ has \newline
order type
$\kappa \times (cf|s^\ell_\eta|)]$.  (This is not really used.)
\endroster
\endproclaim
\bigskip
  
\demo{Proof}  Straightforward and in the most important case see 3.7's
proof.
\enddemo
\bigskip

\remark{Remark}  In the proof we get that each $s^\ell_\nu$ has order
type $\omega$. \newline
\endremark
\bigskip

\proclaim{3.7 Claim}  Suppose in Claim 3.6 that $\kappa = \aleph_0$.
\underbar{Then} we can add
\medskip
\roster
\item "{(i)}"  for $\eta \in S_i,B_\eta$ has the structure of a tree with
$\omega$ levels (e.g., is a family of finite sequences, closed under initial
segments except that $\langle \rangle \notin B_\eta$), and $\eta
\triangleleft \nu
\in S_f$ implies $s^\ell_\eta =  \{ a^\ell_{\eta,m}:m < \omega \}$ is
a branch (of order type $\le \omega$) (a branch is a maximal
linearly ordered subset), and for $m < \ell$, and $k < \omega$, the
$k$'th element of $s^m_\nu$, together with $\nu \restriction \ell$ determines
the $k$-th element of $s^\ell_\nu$.  Also if $\ell < m < n = \ell g(\eta),
\eta \in S_f,\text{cf}[\eta(\ell)] = \lambda(\eta \restriction m) =
\aleph_0$ then $\langle \text{Min}\{ \xi:\text{in } s^\ell_\eta \cap
B_{(\eta \restriction \ell) \char 94 \langle \xi \rangle}$ there are at least
$k$ elements$\}:k < \omega \rangle$ is strictly increasing with limit
$\eta(\ell)$.
\endroster
\endproclaim
\bigskip

\demo{Proof of 3.7}  Without loss of generality let $P$ 
exemplify $PT(\lambda,\kappa)$ fails,
so there are $S$ \, (a $\lambda$-set) and $F,F^0,F^1,F^2$ as in claim 3.3.
As we can shrink $S$, we can assume that it satisfies the conclusion of
3.2(6).
Without loss of generality $\eta \in S_f \Rightarrow lg(\eta) = n$.  
Choose $C_\delta$, $\zeta(\delta,i)$ as 
required in clause (e) (for subclauses $(c), (\alpha), (\beta)$ totally
straight and for subclause $(c)(\gamma)$ we use clause $(e)$ of 3.3).
For $\eta \in S_i$, $\alpha < \lambda(\eta,S)$, we let
$D_{\eta \char 94 \langle \alpha \rangle} =: \cup \{ F^2(\eta \char 94 
< \beta >:
\beta < \alpha,\eta \char 94 < \beta > \in S \}$ so $\langle D_\eta:
\eta \in S_c \rangle$ is a disjoint $\lambda$-system, without loss of
generality disjont to $S$. \newline
For $\eta \in S_f$ and $\ell = 0,\dotsc,n-1$, we define $t^\ell_\eta =:
D_{\eta \restriction(\ell + 1)} \cap \cup \{A:A \in F(\eta)\}$.
\newline
For $\eta \in S_i$ and $\alpha \le \lambda(\eta,S)$ we let
$$
\align
B_{\eta \char 94 \langle \alpha \rangle} =
\biggl\{ \rho:&\rho \text{ is a finite sequence, of length } \ge 3 +
(n - \ell g(\eta)), \\
  &\text{Rang } \rho \subseteq D_{\eta \char 94 \langle \alpha
\rangle} \cup \alpha \cup \{ \eta \} \text{ but Rang}(\rho) \nsubseteq
 \alpha \biggr\}.
\endalign
$$

\noindent
Let
$$
\align
R = \biggl\{ (\ell,m,\eta):&\eta \in S_i, lg(\eta) = m,
\ell \le lg(\eta) \\
 &\text{ and } \lambda(\eta,S) = \text{cf}[\eta(\ell)] > \kappa \biggr\}.
\endalign
$$

\noindent
For $(\ell,m,\eta) \in R$ clearly $\langle \cup \{ t^\ell_\nu:\eta
\triangleleft \nu \in S_f$ and $\nu(m) < \alpha \}:\alpha < 
\lambda(\eta,S) \rangle$ is an
increasing continuous sequence of subsets of $B$ which may have cardinality
$> \lambda(\eta,S)$, each of cardinality
$< \lambda(\eta,S)$.  But $\langle B_{(\eta \restriction \ell) \char 94
\langle \zeta (\eta(\ell),i) \rangle}:i < \lambda(\eta,S) \rangle$ is an
increasing continuous sequence of sets with union
$B_{\eta \restriction(\ell + 1)}$ (remember $\langle
\zeta(\eta(\ell),i):i < \lambda(\eta,S) \rangle$ is an
increasing continuous sequence
of ordinals with limit $\eta(\ell)$ which has cofinality $\lambda(\eta,S)$).
Hence

$$
\align
E_{\eta,\ell} =:
\biggl\{ &i < \lambda(\eta,S):i \text{ is a limit ordinal such that } \\
  &\cup \{ s^\ell_\nu:\eta \triangleleft \nu
\in S_f \} \cap
B_{\eta \restriction (\ell) \char 94 \langle \zeta(\eta(\ell),i) \rangle} \\
  &= \cup \{ s^\ell_\nu:\eta \triangleleft
\nu \in S_f \text{ and } \nu(m) < i \biggr\}
\endalign
$$
\medskip

\noindent
is a club
of $\lambda(\eta,S)$, so let $\varepsilon_{\eta,\ell}:
\lambda(\eta,S) \rightarrow \lambda(\eta,S)$ be a strictly increasing
continuous function with range $E_{\eta,\ell}$.
\newline
It is clear that $\langle B_\eta:\eta \in S_c \rangle$ is a disjoint
$\lambda$-system (note $|B_{\eta \char 94 \langle i \rangle}| < \lambda
(\eta,S)$
as $\lambda(\eta,S)$ is uncountable).  Let $t^\ell_\eta = \{ a(\eta,\ell,i):
i < \omega \}$ (possibly with repetitions).

We define $s^\ell_\eta$ by cases:
\medskip
\roster
\item "{$(\alpha)$}"  if there is $m$ such that $\ell < m < lg(\eta)$,
$(\ell,m,\eta \restriction m) \in R$ and $\lambda(\eta \restriction m,S) >
\aleph_0$ (there is at most one such $m$, and then \newline
$0 \le \ell < m,cf(\eta(\ell)) = \lambda(\eta \restriction m,S) > \aleph_0)$
we let \newline
$\rho^\ell_\eta =: \langle \zeta(n(\ell),\varepsilon_\eta,\ell(\eta(n))
 + 1),\ell,m \rangle \char 94 (\eta \restriction [\ell + 1,n))$,
\newline
$t^\ell_\eta =: \{ \rho^\ell_\eta \char 94 \langle a(\eta,\ell,j):
j \le m \rangle:m < \omega$ and $m > 0\}$
\item "{$(\beta)$}"  $\rho^\ell_\eta = \langle 0,\ell,n \rangle \char 94
\eta \restriction [\ell + 1,n)$, if $cf(\eta(\ell)) \le \kappa$ we let
$$
\align
s^\ell_\eta = \biggl\{ \rho^\ell_\eta \char 94 \langle y_0,\dotsc,y_{2m-1}
\rangle:&m < \omega,m > 0, \text{ for each } k < m, \\
  &y_{2k} = \text{ min} \{ \zeta \in C_{\eta(\ell)}:a(\eta,\ell,0),\dotsc, \\
  &a(\eta,\ell,k) \in B_{\eta \restriction \ell
 \char 94 \langle \zeta \rangle} \}
  \text{and }y_{2k+1} = a(\eta,\ell,k) \biggr\}
\endalign
$$
\endroster

Note that by clause (f) of 3.3, exactly one of those cases occurs.

Now $\langle B_\eta:\eta \in S_c \rangle$, $s^\ell_\eta$ (for $\eta \in S_f,
\ell < \ell g(\eta))$ are as required in 3.6.  The least trivial is (c).  
Suppose
$I \subseteq S_f,|I| < \lambda$, so $\{ \dsize \bigcup_{\ell < n} 
t^\ell_\eta:\eta
\in I \}$ has a transversal, so there is a one-to-one function $g$, 
$\text{ Dom } g = I$ and $g(\eta) \in \dsize \bigcup_\ell t^\ell_\eta$.  
Let $g(\eta) = a(\eta,
h(\eta),g(\eta))$.  Now we define a function $g^*:\text{Dom } g^* = I$,
$g^*(\eta) = \rho^\ell_\eta \char 94 \langle a(\eta,h(\eta),i):0 \le i \le
g(\eta) \rangle$.  Clearly $g^*$ is one-to-one, $g^*(\eta) \in \dsize \bigcup
_{\ell < n} s^\ell_\eta$.

Let for $\eta \in S_i$, $<_\eta$ be a well ordering of
$\{ \eta \} \cup D_{\eta \char 94 < \lambda(\eta,S)>}$ of order type
$\lambda(\eta,S)$ such that $\eta$ is first, and each $\{ \eta \} \cup
D_{\eta \char 94 \langle \alpha \rangle}$ is an initial segment defined by
$\alpha$.  Now $<^*_\eta$ will be $\rho_1 <^*_\eta \rho_2$ iff
$\langle \text{max}_{<_\eta} \text{ Rang }\rho_1 \rangle \char 94 \rho_1
<_{lx} \langle \text{max}_{< \eta} \text{ Rang }\rho_2 \rangle \char 94 \rho_2
<_{lx}$ is lexicographically according to $<_\eta$. \newline
It is also obvious that (i) holds, except possibly the last phrase; but the
correction needed is small so we finish. \hfill$\square_{3.7}$
\enddemo
\bigskip

\proclaim{3.8 Claim}  Suppose $\langle B_\eta:\eta \in S_c \rangle,s^\ell
_\eta(\eta \in S_f,\ell < \ell(\eta))$ are as in Claims 3.6, 3.7; we can
omit 3.6(h)).

Then for any $\rho \in S_i,m = \ell(\rho)$, and $I \subseteq \{ \eta \in
S_f:\rho \le \eta\}$ the following are equivalent:
\medskip
\roster
\item "{$(A)_{\rho,I}$}"  The family $\{ \dsize \bigcup_{\ell \ge m}
s^\ell_\eta:\eta \in I\}$ has a transversal.
\item "{$(B)_{\rho,I}$}"  There are a well ordering $<^*$ of $I$ and
$\{ u_\eta:\eta \in I\}$ such that:
{\roster
\itemitem{ (i) }  for $\eta <^* \nu$ (both in $I$), $u_\nu \cap (\dsize
\bigcup_{\ell \le m} s^\ell_\eta) = \emptyset$.
\itemitem{ (ii) }  For every $\eta \in I$ for some $\ell,m \le \ell <
\ell(\eta),u_\eta$ is an end-segment of $s^\ell_\eta$.
\itemitem{ (iii) }  If $\xi < \text{ Min}\{ \eta(m):\eta \in I\}$ is given,
we can demand that each $u_\eta(\eta \in I)$ is disjoint to
$B_{\rho \char 94 (\xi)}$.
\endroster}
\item "{$(C)_{\rho,I}$}"  There is no $\lambda(\rho,S)$-set $S^*$ such that
$\eta \in S^*_f \Rightarrow \rho \char 94 \eta \in I$.
\item "{$(D)_{\rho,I}$}"  Suppose $\xi < \text{ Min}\{\eta(m):\eta \in I\}$,
there are $u_\eta(\eta \in I)$ where
{\roster
\itemitem{ (i) }  the $u_\eta$ are pairwise disjoint
\itemitem{ (ii) }  $u_\eta$ is an end segment of some $s^\ell_\eta m \le \ell
< \ell(\eta)$
\itemitem{ (iii) }  $u_\eta$ is disjoint to $B_{\rho \char 94 \langle \xi
\rangle}$.
\endroster}
\endroster
\endproclaim
\newpage

\bigskip
REFERENCES
\bigskip

\bibliographystyle{literal-unsrt}
\bibliography{lista,listb,listx}

\enddocument
\bye